\documentclass[12pt]{article}

\usepackage{amsmath}
\usepackage{amssymb}
\usepackage{amscd}
\usepackage{bbm}

\newcommand{\proclaim}[2]{\noindent\textbf{#1.} {\itshape #2}\medskip}
\newcommand{\qed}{\ensuremath{\Box}}

\def\Z{\mathbbm{Z}}
\def\R{\mathbbm{R}}

\def\S{\mathbbm{S}}

\def\D{\mathbbm{D}}
\def\C{\mathbbm{C}}

\def\B{\mathbbm{B}}
\def\a{\alpha}

\def\e{\varepsilon}
\def\d{\delta}

\def\nn {\vskip 0.3cm \noindent }
\def\n {\noindent}

\def\v {\vskip.1cm}
\def\vv {\vskip.2cm}

\makeindex
                              
\begin{document}

\title {\bf  \Large {PROPORTIONALITY OF  INDICES OF  1-FORMS ON  SINGULAR VARIETIES }}

%%\endtitle
\author{\large  J.-P. Brasselet, J. Seade
\thanks{partially supported by CONACYT and DGAPA-UNAM, Mexico}  
 \ and T. Suwa 
\thanks{partially supported by JSPS} } 

%%\endauthor
%%\vskip.3cm
\date {9th May, 2005 }
%\affil
%\linebreak E-mail:  \endaffil

%\date{10 November 2004}

\setcounter{section}{-1}

%%\NoRunningHeads

\pagestyle{headings}

%\pagenumbering{roman}

\maketitle
\bibliographystyle{plain}

\section{Introduction}

M.-H. Schwartz in \cite {Sch1, Sch2} introduced the technique of radial extension of stratified vector fields and frames on singular varieties, and used this to construct cocycles representing  classes in the cohomology
$H^*(M, M\setminus V)$, where $V$ is a singular variety embedded in a complex manifold $M$; these 
are now called {\it the Schwartz classes} of $V$. A basic property of  radial extension is that the index of the vector fields (or frames)  constructed in this way is the same when measured in the strata or in the ambient space; this is called the Schwartz index of the vector field (or frame).
MacPherson in \cite{MP} introduced the notion of the local Euler obstruction, an invariant defined at each point of a singular variety using an index of an appropriate radial 1-form, and used this (among other things) to construct the homology Chern classes of singular varieties. Brasselet and Schwartz in \cite {BS} proved that the Alexander isomorphism
$H^*(M, M\setminus V) \cong H_*(V)$ carries the Schwartz classes into the MacPherson classes; a key ingredient for this proof is their {\it proportionality theorem} relating the Schwartz index and the local Euler obstruction.

These were the first indices of vector fields and 1-forms  in the literature.
 Later in \cite {GSV} was introduced another index for vector fields on isolated hypersurface singularities, and this definition was extended in \cite{SS} to vector fields on complete intersection germs. This is known as the GSV-index and one of its main properties is that it is invariant under perturbations of both, the vector field and the functions that define the singular variety. The definition of this index was recently extended in \cite {BSS1} for vector fields with isolated singularities on hypersurface germs with non-isolated singularities, and it was proved that this index satisfies a proportionality property analogous to the one proved in 
\cite{BS} for the Schwartz index and the local Euler obstruction, the proportionality factor being now the Euler-Poincar\'e characteristic of a local Milnor fiber.

In \cite{EG1} Ebeling and Gusein-Zade observed that when dealing with singular varieties, 1-forms have certain advantages over vector fields, as for instance the fact that for a vector field on the ambient space the condition of being tangent to a (stratified) singular variety is very stringent, while every 1-form on the ambient space defines, by restriction, one on the singular variety.
They adapted the definition of the GSV-index to 1-forms on complete intersection germs with isolated singularities, and proved a very nice formula for it in the case when the form is holomorphic, generalizing the well-known formula of L\^e-Greuel for the Milnor number of a function. 

This article is about 1-forms on complex analytic varieties and it is particularly relevant when the variety has  non-isolated singularities. We show in section 2 how the radial extension technique of M.-H. Schwartz can be adapted to 1-forms,  allowing us to define {\it the Schwartz index} of 1-forms with isolated singularities
on singular varieties. Then we see (section 3) how MacPherson's local Euler obstruction, adapted to 1-forms in general, relates to the Schwartz index, thus obtaining a
proportionality theorem for these indices 
 analogous to the one in \cite{BS} for vector fields. 
  We then extend  (in section 4) the definition of the GSV-index to 1-forms with isolated singularities on (local) complete intersections with non-isolated singularities that satisfy the Thom $a_f$-condition (which is always satisfied if the variety is a hypersurface), and we prove the corresponding proportionality theorem for this index. When the form is the differential of a holomorphic function $h$, this index measures the number of critical points of a generic perturbation of $h$  on a local Milnor fiber, so it is analogous to invariants studied by  a number of authors  (see for instance \cite {Go, IS,STV}).  Section 1 is a review of well-known facts about real and complex valued 1-forms.

The radial extension of 1-forms  can be made global on compact varieties, and it can also be made 
 for frames of differential 1-forms. One gets in this way 
the dual Schwartz classes of singular varieties, which equal the usual ones 
 up to sign. One also has the dual Chern-Mather classes of $V$, already envisaged in \cite{Sa}, and the proportionality formula 3.3 can be used as in \cite{BS} to express  the dual Chern-Mather classes as ``weighted"  dual Schwartz classes, the weights been given by the local Euler obstruction. Similarly,
 in analogy with Theorem 1.1 in \cite {BSS1},  
the corresponding GSV-index and the proportionality Theorem 4.4 extend to frames and can be used to express the dual Fulton-Johnson classes of singular hypersurfaces embedded with trivial normal bundle in compact complex manifolds,
as ``weighted"  dual Schwartz classes, the weights been now given by the Euler-Poincar\'e characteristic of the local Milnor fiber.

 This work was done while the second and third named authors were visiting the 
 ``Institut de Math\'ematiques de Luminy", France; they acknowledge the support  of the CNRS,
 France and  the ``Universit\'e de  la M\'editerran\'ee". 
 
 The authors thank J. Sch\"urmann for his comments and suggestions on the first version of the paper. In particular, he gave us an 
 alternative proof of Theorem 3.3 in the case of the differential form associated to a Morse function, 
 using stratified Morse theory and the micro-local index formula in  \cite {Schu2}.

\section{Some basic facts about 1-forms}

In this section we study some basic facts about the geometry of 1-forms and  the interplay between real and complex valued 1-forms on (almost) complex manifolds, which plays an important role in the sequel.
The material here is all contained in the literature; we include it for completness
and to set up our notation with no possible ambiguities.
  We give precise references when appropriate.

Let $M$ be an almost complex manifold of real dimension $2m>0$. Let $TM$ be its complex tangent bundle.  We denote by $T^*M$ the cotangent bundle of $M$, dual of $TM$; each fiber $(T^*M)_x$ consists of the $\C$-linear maps $TM_x \to \C$. Similarly, we denote by $T_\R M$ the underlying real tangent bundle of $M$; it is a real vector bundle of fiber dimension $2m$, endowed with a canonical orientation. Its dual $T_\R^*M$ has as fiber the $\R$-linear maps $(T_{\R}M)_x \to \R$.

\nn
{\bf  1.1 Definition.} Let $A$ be a subset of $M$. By a real (valued) 1-form $\eta$  on $A$ we mean the restriction to $A$ of a continuous section of the bundle $T_\R^*M$, i.e., for each $x \in A$, 
$\eta_x$ is an $\R$- linear map  $(T_{\R}M)_x \to \R$. We usually drop the word ``valued" here and speak only of real 1-forms on $A$. 
Similarly, a complex  1-form $\omega$ on  $A$ means the restriction to $A$ of a continuous section of the bundle $T^*M$, i.e., for each $x \in A$, 
$\omega_x$ is a $\C$-linear map  $(TM)_x  \to \C$.

\v
Notice that the kernel of a real form $\eta$ at a point $x$ is either the whole fiber $(T_\R M)_x$ or a real 
hyperplane in it. In the first case we say that $x$ is a singular point (or zero) of $\eta$. In the second case
the kernel $ker\, \eta_x$ splits  $(T_\R M)_x$ in two half spaces $(T_\R M^{\pm})_x$; in one of these the form takes positive values, in the other $\eta(v)$ is negative. 

We recall that a vector field $v$ in $\R^N$ is radial at a point $x_o$ if it is transversal to every sufficiently
small sphere around $x_o$ in $\R^N$.  The duality between real 1-forms and vector fields assigns to each 
tangent vector $\partial / \partial x_i$  the form $dx_i$ (extending it by linearity to all tangent vectors). 
This refines the classical duality that assigns to each  hyperplane in 
  $\R^{N}$ the  line orthogonal to it and  motivates the following definition (c. f. \cite {EG1, EG2}):

\nn
{\bf  1.2 Definition.} A real 1-form $\eta$ on $M$ is {\bf radial} (outwards-pointing) at a point $x_o \in M$ if, locally, it is dual over 
$\R$  to a radial outwards-pointing vector field at $x_o$.   Inwards-pointing radial vector fields are defined similarly.

In other words, $\eta$ is {\bf radial} at a point $x_o$ if it is everywhere positive when evaluated in  some radial vector field at $x_o$.

 \v
 
 Thus, for instance, if for a fixed $x_o \in M$  we let $\rho_{x_o}(x)$ be the function $\Vert x - x_o \Vert^2$ (for some Riemmanian metric),  then its differential is a radial form.

\nn
{\bf  1.3 Remark.}
 The concept of radial forms was introduced in \cite{EG1}. In \cite {EG2} 
 radial forms are defined  using more relaxed conditions than we do here. 
 However this is a concept "imported" from the corresponding notion of radial vector fields, so
  we use definition 1.2.  
 \v
 
 A complex 1-form $\omega$ on $A \subset M$ can be written in terms of its real and imaginary parts:
 $$\omega \,=\, Re\,(\omega) \,+\, i\, Im \,(\omega)\,.$$
 Both $ Re\,(\omega)$ and $ Im \,(\omega)$ are real 1-forms, and the linearity of $\omega$ implies that for each tangent vector one has:
 $$  Im \,(\omega)(v) \,=\, - Re \, (\omega)(iv) \,,$$
 thus 
 \[\omega (v) \,=\, Re\,(\omega) (v) \,-\, i\, Re \,(\omega) (iv)\,. \]
 In other words the form $\omega$ is determined by its real part and one has a 1-to-1 correspondence between real and complex forms, assigning to each complex form its real part, and conversely, to a real 1-form $\eta$ corresponds the complex form $\omega$ defined by:
 \[\omega(v) \,=\, \eta(v) - i \eta(iv)\,. \]

 This statement (noted in \cite {EG2},\cite{GMP}) 
  refines the obvious fact that a complex hyperplane $P$ in $\C^m$, say defined by a linear form $H$,  is the intersection of
the real hyperplanes $\widehat H := \{Re\, H = 0\}$ and $\,i\, \widehat H$. This justifies the following definition:
 
 \nn
{\bf  1.4 Definition.} A complex 1-form $\omega$ is {\bf radial} at a point $x \in M$ if its real part is radial at $x$.

  \v

Recall that the Euler class of an oriented vector bundle 
is the primary obstruction for constructing a non-zero section \cite {St}. 
In the case of the bundle $T^*_\R M$,
this class equals the Euler class $\hbox{Eu}(M)$ of the underlying real tangent bundle $T_\R M$, since they are isomorphic.  Thus, if $M$ is compact  then its Euler class evaluated on the orientation cycle of 
$M$ gives the Euler-Poincar\'e characteristic $\chi(M)$. We can say this in different words: let 
$\eta$ be a real 1-form on $M$ with isolated (hence finitely many) singularities $x_1, \cdots,x_r$. At
each $x_i$ this 1-form defines a map, 
$\S_\e \buildrel{\eta/\Vert \eta \Vert}\over \longrightarrow \S^{2m-1}$, from a small sphere in $M$ around
$x_i$ into the unit sphere in the fiber $(T^*_\R M)_x$. The degree of this map is the 
{\bf Poincar\'e-Hopf} local index of $\eta$ at $x_i$, that we may denote by 
$\hbox{Ind}_{PH}(\eta, x_i)$.
 Then the total index of $\eta$ in $M$ is by definition the sum of its local indices at the $x_i$ and it equals $\chi(M)$. Its Poincar\'e dual class in
$H^{2m}(M)$ is the Euler class of $T^*_\R M \cong T_\R M$. 

\v

More generally, if $M$ is a compact, $C^\infty$ manifold of real dimension $2m$ with 
non-empty boundary $\partial M$ and a complex structure in its tangent bundle, one can speak of real and complex valued 1-forms as above. 
Elementary obstruction theory (see \cite {St}) implies that one can always find real and complex 
1-forms on $M$ with isolated singularities, all contained in the interior of $M$. In fact, if a 
real 1-form
$\eta$ is defined in a neighborhood of $\partial M$ in $M$ and it is non-singular there, then we can always extend it to the interior of $M$ with finitely many singularities, and its total 
 index in $M$  does not depend on the choice of the extension. 
 
  \nn
{\bf  1.5 Definition.} Let $M$ be an almost complex manifold with boundary $\partial M$ and let 
$\omega$ be  a (real or complex) 1-form on $M$, non-singular on a neighborhood of $\partial M$;
let $Re \, \omega$ be its real part if $\omega$ is a complex form, otherwise $Re \, \omega = \omega$ for real forms. 
The form $\omega$ is {\bf radial} at the boundary if for each  vector $v(x) \in TM$, $x \in \partial M$, 
which is normal to the boundary (for some metric), 
pointing outwards of $M$, one has $Re \, \omega(v(x)) > 0$.
 
 \vv
 By the 
 theorem of Poincar\'e-Hopf for manifolds with boundary, if a real 1-form $\eta$ is radial at the boundary and $M$ is compact, then the total index of $\eta$ is $\chi(M)$. 
 
 \v
 
We now make similar considerations for complex 1-forms. We let $M$ be a compact, $C^\infty$ manifold of real dimension $2m$ (with or without boundary $\partial M$), with 
 a complex structure in its tangent bundle $TM$. Let $T^*M$ be as before, the cotangent bundle of $M$, i.e., the bundle of complex valued continuous 1-forms. The top Chern class $c^m(T^*M)$ is the
primary obstruction for constructing a section of this bundle, i.e., if $M$ has empty boundary, then
$c^m(T^*M)$ is the number of points, counted with their local indices, of the zeroes of a section 
$\omega$ of $T^*M$  (i.e., a complex 1-form) with isolated singularities (i.e., points where it vanishes). 
It is well known (see for instance \cite {Mi}) that one has:
$$c^m(T^*M) \,=\, (-1)^m \,c^m(TM)\,.$$
This corresponds to the fact that at each isolated singularity $x_i$  of $\omega$ one has two local indices: one of them is the index of its real part defined as above, $\hbox{Ind}_{PH}(Re \, \omega, x_i)$;  the other is the degree of the map 
 $\S_\e \buildrel{\omega/\Vert \omega \Vert}\over \longrightarrow \S^{2m-1}$, that we denote by
 $\hbox{Ind}_{PH}(\omega, x_i)$. These two indices are related by the equality:
 \[\hbox{Ind}_{PH}(\omega, x_i) \,=\, (-1)^m \, \hbox{Ind}_{PH}(Re \, \omega, x_i)\,, \]
 and the index on the right corresponds to the local Poincar\'e-Hopf index of the vector field defined
by duality near  $x_i$.  For instance, the form $\omega = \sum z_i dz_i$ in $\C^m$ has index $1$ at
$0$, while its real part $\sum(x_i dx_i - y_i dy_i)$ has index $(-1)^m$.
 
 If we take $M$ as above, compact and with possibly non-empty boundary, and $\omega$ is a complex 1-form
 with isolated singularities in the interior of $M$ and radial on the boundary, then (by the previous considerations) the total index of $\omega$ in $M$ is $(-1)^m \, \chi (M)$.  We summarize some of the previous discussion in the following theorem (c.f. \cite {EG1, EG2}):
 
 \nn
 \proclaim{1.6 Theorem}{Let $M$ be a compact, $C^\infty$ manifold of real dimension $2m$ (with or without boundary $\partial M$), with  a complex structure in its tangent bundle $TM$. Let $T^*_\R M$ and  $T^*M$ be as before, the bundles of real and complex valued continuous 1-forms on $M$, respectively. Then:
 
 \v \n {\bf i) } Every real 1-form $\eta$ on $M$ determines a complex 1-form $\omega$ by the formula
 $$ \omega(v) \,=\, \eta(v) - i \eta(iv) \,;$$
so the real part of $\omega$ is $\eta$.

 \v \n {\bf ii) }  The local Poincar\'e-Hopf indices at an isolated  singularity of a complex 1-form and its real part are related by:
 \[
 {\rm{Ind}}_{PH}(\omega, x_i) \,=\, (-1)^m \, {\rm{Ind}}_{PH}(Re \, \omega, x_i)\,.
 \]

 \v \n {\bf iii) } If a real 1-form on $M$ is radial at the boundary $\partial M$, then its total Poincar\'e-Hopf
 index in $M$ is $\chi(M)$.  In particular, a radial real 1-form has local index 1.
 
 \v \n {\bf iv) } If a complex 1-form on $M$ is radial at the boundary $\partial M$, then its total Poincar\'e-Hopf  index in $M$ is $(-1)^m  \chi(M)$.
 }

 \nn
 {\bf 1.7 Remark.} One may consider frames of complex 1-forms on $M$ instead of a single 1-form. This means considering sets of $k$ complex 1-forms,  whose singularities are the points where these forms become linearly dependent over $\C$. By definition (see \cite {St}) the primary obstruction for constructing such a frame is the Chern class $c^{m-k+1}(T^*M)$, so these classes also have an expression similar to 1.6 but using indices of frames of 1-forms.  One always has 
 $c^i(T^*M) = (-1)^i c^i(TM)$.  Thus the Chern classes, and  all the Chern numbers of $M$, can be computed using indices of either  vector fields or 1-forms.

\section{Radial extension and the Schwartz index}

 In the sequel we will be interested in considering forms defined on singular varieties in a complex manifold, so we introduce some standard notation.
 Let $V$ be a reduced, equidimensional complex analytic space of dimension $n$ in a complex manifold $M$ of dimension $m$, 
 endowed with a Whitney stratification $\{V_\a\}$ adapted to $V$, i.e., $V$ is union of strata. 

The following definition is an immediate extension for 1-forms of the corresponding (standard)  
definition for  functions on stratified spaces 
in terms of its differential (c.f. \cite {EG2, GMP, Le1}).

\nn
{\bf 2.1 Definition. } Let $\omega$ be a (real or complex) 1-form on $V$,
 i.e., a continuous section of either $T^*_\R M \vert_V$ or $T^*M \vert_V$.   
A singularity of $\omega$ with respect to the Whitney stratification $\{V_\a\}$ means a point  $x$ where
the kernel of $\omega$  contains the tangent space of the corresponding stratum.

This means that the pull back of the form to $V_\a$ vanishes at $x$. 

In section 1 we introduced the notion of radial forms, which is dual to
the "radiality" for vector fields.  We now extend this notion relaxing the condition of radiality in the
directions tangent to the strata.  From now on, unless it is otherwise stated explicitely, 
 by a singularity of a 1-form on $V$ we mean a singularity in the stratified sense, i.e., in the sense of 2.1.

\nn
{\bf 2.2 Definition. } Let $\omega$ be a (real or complex) 1-form on $V$. The form is {\bf normally radial}
at a point $x_o \in V_\a \subset V$ if it is radial when restricted to vectors which are not tangent to the stratum $V_\a$ that contains $x_o$. In other words,  for every vector $v(x)$ tangent to $M$ at a point $x \notin V_\a$, $x$ sufficiently close to $x_o$  and $v(x)$ pointing outwards a tubular neighborhood of the stratum $V_\a$, one has  $Re \, \omega(v) > 0$ (or $Re \, \omega(v) < 0$ for all such vectors; 
if $\omega$  is real then it  equals  $Re \, \omega$).

\vv
Obviously a radial 1-form is also normally radial, since it is radial in all directions.  

For each point $x$ in a stratum $V_\alpha$, one has a neighborhood $U_x$ of $x$ in $M$ 
which is diffeomorphic to the product
$U_\a \times \D_\a,$ 
where $U_\a = U_x \cap V_\a$ and $\D_\a$ is a small disc in $M$ transversal to $V_\a$. 
Let $\pi$ be the projection $\pi : U_x \to U_\a$ and $p$  the projection $p: U_x \to \D_\a$.
One has an isomorphism:
\[\;T^* U_x  \,\cong 
\pi^* T^*{U_\a} \oplus p^* T^*\D_\a\,.\]

That a (real or complex) 1-form $\omega$ be normally radial at $x$  means that up to a local change of coordinates in $M$,  $\omega$ is the direct sum of the pull back of a (real or complex) form on $U_\a$, i.e., a section
of the (real or complex) cotangent bundle $T^*U_\a$, and a  section
of the (real or complex) cotangent bundle $T^*\D_\a$ which is a radial form in the disc. 
\v

It is possible to make for 1-forms 
 the classical construction of {\bf radial extension} introduced by M.-H. Schwartz in \cite{Sch1, Sch2} for stratified vector fields and frames. Locally, the construction can be described  as follows. 
 We consider first real 1-forms.
Let $\eta$ be a 1-form on  $U_\a$, denote by $\widehat \eta$ its pull back to a section of $\pi^* T^*_{\R}{U_\a}$. This corresponds to the  {\bf parallel extension} of stratified vector fields
done by Schwartz. Now look at the function $\rho$ 
given by the square of the distance to the origin in $\D_\a$.  The form $p^* d\rho$ on $U_x$ 
vanishes on $ U_\a $ and away from $U_\a $ its
 kernel is transversal to the strata of $V$ by Whitney conditions.
  
 The sum $\eta' = \widehat \eta + p^*d\rho$ defines a normally radial 1-form on $U_x$ which coincides with $\eta$ on $U_\a $; away from $U_\a $ its kernel is transversal to the strata of $V$. Thus,
if $\eta$ is non-singular at $x$, then $\eta'$ is non-singular everywhere on
 $U_x$. If $\eta$ has an isolated singularity at $x \in V_\a$, then $\eta'$ also has an isolated singularity there.  In particular, if the dimension of the stratum $V_\a$ is zero then $\eta'$ is 
 a radial form in the sense of section 1. 
 
 Following the terminology of \cite{Sch1, Sch2} we say that the form $\eta'$ is obtained from
 $\eta$ by {\bf radial extension}. 
 
Since the index
in $M$ of a normally radial form is its index in the stratum times the index of a radial form in
the disc $\D_\a$, we obtain the following important property of forms constructed by radial extension. 

 \nn
 \proclaim{2.3 Proposition}{Let $\eta$ be a real 1-form on the stratum $V_\a$ with an isolated singularity
 at a point $x$ with local Poincar\'e-Hopf index  ${\rm{Ind}}_{PH}(\eta, V_\a; x)$. Let $\eta'$ the 1-form on
 a neighborhood of $x$ in $M$ obtained by radial extension. Then the 
 index of $\eta$ in the  stratum equals the index of $\eta'$   in $M$:
 $${\rm{Ind}}_{PH}(\eta, V_\a; x)\,=\, {\rm{Ind}}_{PH}(\eta', M; x)\,. $$}

 \n{\bf  2.4 Definition.} The {\bf Schwartz index} of the continuous real 1-form $\eta$ at  an isolated
singularity $x \in V_\a \subset V$, denoted $\hbox{Ind}_{Sch}(\eta, V; x)$, 
 is the Poincar\'e-Hopf index  of the 1-form $\eta'$ obtained from $\eta$ by radial extension; or equivalently, if  the stratum of $x$ has dimension more than 0, 
  $\hbox{Ind}_{Sch}(\eta, V; x)$ is the Poincar\'e-Hopf index  at $x$ of  
 $\eta$ in the stratum $V_\a$.
 
 \vv
 If $x$ is an isolated singularity of $V$ then every 1-form on $V$ must be singular at $x$ since its kernel  
 contains the ``tangent space" of the stratum. In this case the 
 index of the form in the stratum is defined to be 1, and this is consistent with the previous definition since in this case the radial extension of $\eta$ is actually radial at $x$, so it has index 1 in the ambient space.

The previous process is easily adapted to give radial extension for 
complex 1-forms. Let $\omega$ be such a form  on $V_\a$;  let $\eta$ be its real part. We extend $\eta$ as above, by radial extension, to obtain a real 1-form $\eta'$ which is normally radial at $x$. Then we use statement i) in Theorem 1.6 above to obtain a complex 1-form   $\omega'$ on $U_x $ 
that extends $\omega$ and is also normally radial at $x$. If we prefer, we can make this process in a different but equivalent way: first make a parallel extension of
$\omega$ to $U_x $ as above, using the projection $\pi$; denote
by $\widehat \omega$ this complex 1-form. Now use 1.6.i)  to define a complex 1-form 
$\widehat { d\rho}$ on   $ U_x $  whose real part is $d\rho$, and take the direct sum of 
$\widehat \omega$ and 
$\widehat  {d\rho}$ at each point to obtain the extension $\omega'$. We say that 
$\omega'$ is obtained from $\omega$ by {\bf radial extension}. 

We have the equivalent of Proposition 2.3 for complex forms, modified with the appropriate signs:
 $$(-1)^s \, \hbox{Ind}_{PH}(\omega, V_\a; x)\,=\, (-1)^m \, \hbox{Ind}_{PH}(\omega', M; x)\,, $$
where $2s$ is the real dimension of $V_\a$ and $2m$ that of $M$.

\nn {\bf  2.5 Definition.} The {\bf Schwartz index} of the continuous complex 1-form $\omega$ at  an isolated
singularity $x \in V_\a \subset V$, denoted $\hbox{Ind}_{Sch}(\omega, V; x)$,  is $(-1)^n$-times the index of its real part:
$$\hbox{Ind}_{Sch}(\omega, V; x) \,=\, (-1)^n \hbox{Ind}_{Sch}(Re\, \omega, V; x)\;.$$

\section[Local Euler obstruction]
{Local Euler obstruction and the Proportionality Theorem}

We are now concerned only with a local situation, so we take $V$ to be embedded in an
open ball $\B \subset \C^m$ centered at the origin $0$.  On the regular part of $V$ one has the  map 
$\sigma : V_{reg} \to G_{n,m}$
into the Grassmannian of $n (= \text{dim}\,V)$-planes in $\C^m$, that assigns to each point the corresponding tangent space of $V_{reg}$. The 
Nash blow up $\widetilde V \buildrel {\nu}\over \to V$ of $V$ 
is by definition the closure in $\B \times G_{n,m}$ of the
graph of the  map $\sigma$. 
One also has the Nash bundle $\widetilde T \buildrel {p}\over \to \widetilde V $, 
restriction to $\widetilde V $ of the tautological bundle over $\B \times G_{n,m}$.

The corresponding dual bundles of complex and real 1-forms  are denoted by
$\widetilde T^* \buildrel {p}\over \to \widetilde V $ and 
$\widetilde T^*_\R \buildrel {p}\over \to \widetilde V $,  respectively. Observe that a point in ${\widetilde T^*}$ is a triple $(x,P,\omega)$ where $x$ is in $V$, $P$ is an $n$-plane in the tangent space $T_x \B$ which is limit of a sequence $\{(T V_{reg})_{x_i}\}$, where the $x_i$ are points in the regular part of $V$ converging to $x$, and $\omega$ is a $\C$-linear map $P \to \C$.  (Similarly for $\widetilde T^*_\R$.) 

 %\nn
% \proclaim{3.1 Proposition}{Given 1-form .... (canonical lift)....
% }
 
 Let us denote by $\rho$ the function given by the square of the distance to $0$.
We recall that MacPherson in \cite {MP} observed that  the Whitney condition (a) implies that the pull-back of the differential $d\rho$ defines a never-zero section $\widetilde  {d\rho}$ of
 $\widetilde T^*_\R$ over $\nu^{-1}(\S_\e \cap V) \subset \widetilde V$, where $\S_\e$ is the boundary of a small ball $\B_\e$ in $\B$ centered at $0$. The obstruction for extending $\widetilde  {d\rho}$ as a
 never-zero section of $\widetilde T^*_\R$ over $\nu^{-1}(\B_\e \cap V) \subset \widetilde V$ is a 
 cohomology class in 
 $H^{2n}(\nu^{-1}(\B_\e \cap V), \nu^{-1}(\S_\e \cap V); \Z)$, and 
 MacPherson defined {\bf the local Euler obstruction} $\hbox{Eu}_V(0)$ of $V$ at
 $0$  to be the integer  obtained by evaluating  this class 
on the orientation cycle
 $[\nu^{-1}(\B_\e \cap V), \nu^{-1}(\S_\e \cap V)]$.

More generally, given a section $\eta$ of $T^*_{\Bbb R}\Bbb B|_A$, $A\subset V$, there is a canonical way of constructing a section $\tilde\eta$ of $\widetilde T^*_\R|_{\tilde A}$, 
$\tilde A=\nu^{-1}A$, which is described in the following. The same construction 
works for complex forms. First, taking the pull-back $\nu^*\eta$, we get a section of 
$\nu^*T^*_{\Bbb R}\Bbb B|_V$. Then $\tilde\eta$ is obtained by projecting  $\nu^*\eta$
to a section of $\tilde T^*_\R$ by the canonical bundle homomorphism
$$
\nu^*T^*_{\Bbb R}\Bbb B|_V\longrightarrow\tilde T^*_\R.
$$
Thus the value of $\tilde\eta$ at a point $(x,P)$ is simply the restriction of the linear
map $\eta(x): (T_{\Bbb R}\Bbb B)_x\to\Bbb R$ to $P$. 
We call $\tilde\eta$ the {\bf canonical lifting} of $\eta$.

%More generally, a  (real or complex) 1-form $\eta$ defined on a subset $A \subset V$ has a canonical lifting to a (real or complex, respectively) section 
%$\widetilde \eta$ of  $\widetilde T^*_\R$  (respectively $\widetilde T^*$) over $\nu^{-1}(A)$: a fiber of this bundle is by definition
%a $n$-plane $P$ in $\B$ which is limit of a sequence of planes tangent to the regular part of $V$; 
 %for each point $\tilde x = (x,P)$ in 
%$\nu^{-1}(A)$, one defines $\widetilde \eta(\tilde x) = (x,P,\eta)$. 
By the Whitney condition (a), if 
$a \in V_\alpha$ is the limit point of the sequence $\{ x_i \}\in V_{\rm reg}$ such that $P =
 \lim (TV_{\rm reg})_{x_i} $ and if
the kernel of $\eta$ is transversal to $V_\alpha$, then the linear form $\widetilde \eta$ will be non-vanishing on $P$. Thus, if $\eta$
has an isolated singularity
at the point $0 \in V$ (in the stratified sense), then we have a never-zero section $\widetilde \eta$ 
of the dual Nash bundle $\widetilde T^*_\R$ over $\nu^{-1}(\S_\e \cap V) \subset \widetilde V$. Let 
$o(\eta) \in H^{2n}(\nu^{-1}(\B_\e \cap V), \nu^{-1}(\S_\e \cap V); \Z)$ be the 
cohomology class of the obstruction cycle 
to extend this to a section of  $\widetilde T^*_\R$ over $\nu^{-1}(\B_\e \cap V)$. Then define (c.f. \cite {BMPS, EG2}):

\nn
{\bf 3.1 Definition.} The {\bf local Euler obstruction} of the real differential form $\eta$ at an isolated singularity is the integer 
$\hbox{Eu}_{V}(\eta,0)$ 
obtained by evaluating the obstruction cohomology class $o(\eta)$ on the orientation cycle
$[\nu^{-1}(\B_\e \cap V), \nu^{-1}(\S_\e \cap V)]$.
\nn

The local Euler obstruction $\hbox{Eu}_{V}(0)$ 
of MacPherson   corresponds to taking  the differential of the squared function distance to $0$. 

In the complex case, one can perform the same construction, using the 
corresponding complex bundles. If $\omega$ is a complex differential form,
section of $T^*\Bbb B|_A$ with an isolated singularity, one can define  the local Euler obstruction 
$\hbox{Eu}_{V}(\omega,0)$.
 Notice that it is equal to  that of its real part up to sign:
\[\hbox{Eu}_{V}(\omega,0) \,=\, (-1)^n \hbox{Eu}_{V}(Re\, \omega, 0) \,. \tag{3.2}\]
This is an immediate consequence of the relation between the Chern classes of a complex vector bundle and those of its dual (see for instance \cite {Mi}).

We note that the idea to consider the (complex) dual Nash bundle was already present in
\cite{Sa},  where Sabbah  introduces a local  Euler obstruction ${\rm E\check u}_V (0)$ that satisfies
${\rm E\check u}_V (0) = (-1)^{n}{\rm Eu}_V(0)$. See also Sch\"urmann \cite{Schu1}, sec. 5.2.

\v
Just as for vector fields (see \cite {BS}), one has in this situation the following:

\nn
\proclaim{3.3 Theorem} {Let  $ V_\a \subset V$ be the stratum containing $0$, 
${\rm{Eu}}_V(0)$ the local Euler obstruction of $V$ at $0$ and 
$\omega$ a (real or complex) 1-form on $V_\a$  with an isolated singularity at $0$. 
Then the  local Euler obstruction of the
  radial extension $\omega'$ of $\omega$ and the Schwartz index of $\omega$ at $0$ are related by the 
following   proportionality formula:
  
  \[{\rm{Eu}}_{V}(\omega',0) \,=\, {\rm{Eu}}_{V} (0)\cdot {\rm{Ind}}_{Sch}(\omega, V; 0) \,.\]}

\nn
{\bf Proof} 
By 3.2 and Theorem 1.6 above, it is enough to prove 3.3 for either real or complex 1-forms, each case implying the other. We prove it for real forms.

Let $\eta$ and   $\eta'$  be as above. Also, let $\eta_{rad}$ 
denote a radial form at $0$. 
  
By construction and definition, we have
\[
\hbox{Ind}_{PH}(\eta,V_\alpha; 0)=\hbox{Ind}_{PH}(\eta',\Bbb B;0)=\hbox{Ind}_{Sch}(\eta,V;0).\tag{3.4}
\]

By definition of $\hbox{Ind}_{PH}(\eta',\Bbb B;0)$, there is a homotopy
$$
\Psi:[0,1]\times\Bbb S_\e\longrightarrow T^*_\R\Bbb B|_{\Bbb S_\e}
$$
such that its image satisfies:
\[
\partial\hbox{Im}\Psi=\eta'(\Bbb S_\e)-\hbox{Ind}_{PH}(\eta',\Bbb B;0)
\cdot \eta_{rad}(\Bbb S_\e)\tag{3.5}
\]
as chains in $T^*_{\Bbb R}\Bbb B|_{\Bbb S_\e}$. The restriction of $\Psi$ gives a homotopy
$$
\psi:[0,1]\times(\Bbb S_\e\cap V)\longrightarrow T^*_{\Bbb R}\Bbb B|_{\Bbb S_\e\cap V}
$$
such that (c.f. (3.4))
$$
\partial\hbox{Im}\psi=\eta'(\Bbb S_\e\cap V)-\hbox{Ind}_{Sch}(\eta,V;0)
\cdot\eta_{rad}(\Bbb S_\e\cap V).
$$

Now we can lift $\psi$, $\eta'$ and $\eta_{rad}$ to sections 
$\nu^*\psi$, $\nu^*\eta'$ and $\nu^*\eta_{rad}$ of $\nu^*T^*_{\Bbb R}\Bbb B$ to get
a homotopy
$$
\nu^*\psi:[0,1]\times\nu^{-1}(\Bbb S_\e\cap V)
\longrightarrow\nu^*T^*_{\Bbb R}\Bbb B|_{\nu^{-1}(\Bbb S_\e\cap V)}
$$
and, since $\nu$ is an isomorphism away from the singularity of $V$, we still have
\[
\partial\hbox{Im}\nu^*\psi=\nu^*\eta'(\nu^{-1}(\Bbb S_\e\cap V))-\hbox{Ind}_{Sch}(\eta,V;0)
\cdot\nu^*\eta_{rad}(\nu^{-1}(\Bbb S_\e\cap V))
\]
as chains in $\nu^*T^*_{\Bbb R}\Bbb B|_{\nu^{-1}(\Bbb S_\e\cap V)}$. Recall that we get the canonical liftings
$\tilde\psi$, $\tilde\eta'$ and $\tilde\eta_{rad}$ of $\psi$, $\eta'$ and $\eta_{rad}$
by taking the images of $\nu^*\psi$, $\nu^*\eta'$ and $\nu^*\eta_{rad}$ by the canonical
bundle homomorphism
$\nu^*T^*_{\Bbb R}\Bbb B \longrightarrow \tilde T^*_{\Bbb R}$.
Thus we have
\[
\partial\hbox{Iml}\tilde\psi=\tilde\eta'(\nu^{-1}(\Bbb S_\e\cap V))-\hbox{Ind}_{Sch}(\eta,V;0)
\cdot\tilde\eta_{rad}(\nu^{-1}(\Bbb S_\e\cap V))
\]
as chains in 
$\tilde T^*_{\Bbb R}|_{\nu^{-1}(\Bbb S_\e\cap V)}$.
The forms $\tilde\eta'$ and $\tilde\eta_{rad}$ are non-vanishing on 
$\nu^{-1}(\Bbb S_\e\cap V)$, by the Whitney condition, and by definition of the Euler obstructions, we have 
the theorem.
\qed

\section{The GSV-index}

We recall (\cite {GSV, SS}) that the GSV-index of a vector field $v$ on an isolated complete intersection germ $V$ can be defined to be the Poincar\'e-Hopf index of an extension of $v$ to a Milnor fiber $F$. Similarly, the GSV-index of a 1-form $\omega$ on $V$ can be defined to be the Poincar\'e-Hopf index of the form on $F$, i.e.,  the number of singularities of $\omega$ in $F$ counted with multiplicities \cite {EG1}. When $V$ has 
non-isolated singularities one may not have a Milnor fibration in general, but one does if $V$ has
a Whitney stratification with Thom's $a_f$-condition, $f =(f_1,\cdots,f_k)$ being the functions that define $V$ (c.f. \cite {Le2, LT, BSS1}).

Let $(V,0)$ be a  complete intersection of complex dimension $n$ defined in a ball $\B$ in
$\C^{n+k}$  by functions
$f =(f_1,\cdots,f_k)$, and assume $0$ is a singular point of $V$ (not necessarily an isolated singularity). As before, we endow $\B$ with a Whitney stratification adapted to $V$, and we assume that we can choose
$\{V_\a\}$ so that it satisfies the $a_f$-condition of Thom (see for instance \cite{LT}). In particular one always has such a stratification if $k = 1$, by \cite{Hi}. 

Let $\omega$ be as before, a (real or complex) 1-form on $\B$, and assume  its restriction
to $V$ has an isolated singularity at $0$. This means that the kernel of $\omega(0)$ contains the tangent space of the stratum $V_\a$ containing $0$,
 but everywhere else it is transversal to each stratum
$V_\a \subset V$. Now let $F = F_t$ be a Milnor fiber of $V$, i.e., $F = f^{-1}(t) \cap \B_\e$, where
$\B_\e$ is a sufficiently small ball in $\B$ around $0$ and $t \in \C^k$ is a regular value of $f$ with
$\Vert t \Vert$ sufficiently small with respect to $\e$. Notice that the $a_f$-condition implies that for every sequence $t_n$ of regular values converging to $0$, and for every sequence $\{x_n\}$ of points in the corresponding Milnor fibers converging to a point $x_o \in V$ so that the sequence of tangent spaces
$\{(TF)_{x_n}\}$ has a limit $T$, one has that  $T$ contains the space 
$(TV_\a)_{x_o}$, tangent to the stratum that contains $x_o$. By transversality this implies that 
choosing the regular value $t$ sufficiently close to $0$ we can assure that the kernel of
$\omega$ is transversal to the Milnor fiber at every point in its boundary $\partial F$. Thus its pull-back
to $F$ is a 1-form on this smooth manifold, and it is never-zero on its boundary, thus
$\omega$ has a well defined Poincar\'e-Hopf index in $F$ as in section 1. This index is well-defined and depends only on the restriction of
$\omega$ to $V$ and the topology of the Milnor fiber $F$, which is well-defined once we fix the defining function $f$ (which is assumed to satisfy the $a_f$-condition for some Whitney stratification). 

\nn
{\bf 4.1 Definition} The GSV-index of $\omega$ at $0 \in V$ relative to $f$, $\hbox{Ind}_{GSV}(\omega,0)$, is the 
Poincar\'e-Hopf index of $\omega$  in $F$.

\vv
In other words this index measures the number of points (counted with signs) in which a generic perturbation of $\omega$  is tangent to $F$.
In fact the inclusion $F \buildrel{i}\over \to M$ pulls the form $\omega$ to a section of the
(real or complex, as the case may be) cotangent bundle of $F$, which is never-zero near the 
boundary because $\omega$ has an isolated singularity at $0$ and, by hypothesis,  the map $f$ satisfies the $a_f$-condition of Thom.  If the form $\omega$ is real then 
\[ \hbox{Ind}_{GSV}(\omega,0) \,=\, \hbox{Eu}(F; \omega)[F] \;, \tag{4.2}\]
where $\hbox{Eu}(F; \omega) \in H^{2n}(F, \partial F)$ is the Euler class of the real cotangent bundle 
$T^*_\R F$ relative to the section defined by $\omega$ on the boundary, and $[F]$ is the orientation cycle of the pair $(F,\partial F)$. 
If $\omega$ is a complex form, then one has:
\[ \hbox{Ind}_{GSV}(\omega,0) \,=\, c^n(T^*F; \omega)[F] \;, \tag{4.3.i}\]
where $c^n(T^*F; \omega)$ is the top Chern class of the cotangent bundle of $F$ relative to the form $\omega$ on its boundary. In this case one can, alternatively,  express this index as the  relative Chern class:
\[ \hbox{Ind}_{GSV}(\omega,0) \,=\, c^n(T^*M \vert_F; \Omega)[F] \;, \tag{4.3.ii}\]
where $\Omega$ is the frame of $k+1$ complex 1-forms on the boundary of $F$ given by
\[\Omega \,=\, (\omega, df_1, df_2, \cdots, df_k)\,,\]
since the forms 
$(df_1,\cdots,df_k)$ are linearly independent everywhere on $F$. 
Notice that if the form $\omega$ is holomorphic, then this index 
is necessarily non-negative  because it can be regarded as an intersection number of complex submanifolds.  For every complex 1-form one has:
\[ \hbox{Ind}_{GSV}(\omega,0) \,=\, (-1)^n  \hbox{Ind}_{GSV}(Re\, \omega,0) \,.\]

We remark that if $V$ has an isolated singularity at $0$, this is the index envisaged in \cite{EG1}, i.e.,
the degree of the map from the link $K$ of $V$ into the Stiefel manifold of complex (k+1)-frames
in the dual $\C^{n+k}$ given by the map $(\omega, df_1,\cdots,df_k)$.  Also notice that this index is somehow dual to the index defined in \cite {BSS1} for vector fields, which is related to the top Fulton-Johnson class of singular hypersurfaces.

So, given the (non-isolated) complete intersection singularity $(V,0)$ and a (real or complex) 1-form
$\omega$ on $V$ with an isolated singularity at $0$, one has three different indices:  the Euler obstruction (section 2),  the GSV-index just defined and 
  the index of its pull back to a 1-form on the stratum of $0$. One also has the index of the form in the ambient manifold $M$.
  For forms obtained by radial extension, the index in the stratum equals its index in $M$, and this is by definition the Schwartz index.
   The following proportionality theorem is analogous to the one in \cite {BSS1} for vector fields.

\nn
\proclaim{4.4 Theorem} {Let $\omega$ be a (real or complex) 1-form on the  stratum $V_\a$ of $0$ with an isolated singularity at $0$.  Then the GSV index of its radial extension $\omega'$ 
 is proportional to the Schwartz index,  the proportionality factor being the Euler-Poincar\'e characteristic of the Milnor fiber $F$:
\[{\rm{Ind}}_{GSV}(\omega',0)\,=\, {\chi(F) \cdot \rm{Ind}}_{Sch}(\omega, V; 0) \,.\]
}

\v \n
{\bf Proof.} 
It is enough to prove 4.4 either for complex forms or for real forms, each one implying the other.  The proof is similar to that of 3.3.
  Let $\omega'$ and $\omega_{rad}$ be as in the proof of Theorem 3.3. Then 4.4 is proved by
taking the retriction to $F$ of each section in (3.5) as a differential form, noting that $\hbox{Ind}_{GSV}(\omega_{rad},0)\,=\, \chi(F)$. \qed

\nn
{\bf 4.5 Remark.}  We notice that 4.2 and 3.3 can also be proved using the stability of the index under
perturbations; this  works for vector fields too. More precisely, one can easily show that the Euler obstruction $\hbox{Eu}_V(\omega,x)$ and the GSV-index are stable when we perturb the 1-form (or the vector field) in the stratum and then extend it radially; then the sum of the indices at the singularities of the new 1-form (vector field) give the corresponding index for the original singularity. This implies the proportionality of the indices.

\end{document}